\documentclass[12pt]{article}
\makeatletter
  
  \@addtoreset{equation}{section}
\makeatother

\newtheorem{thm}{Theorem}[section]
\newtheorem{lm}{Lemma}[section]
\newtheorem{rk}{Remark}[section]
\newtheorem{co}{Corollary}[section]
\newtheorem{df}{Definition}[section]

\usepackage{amsmath,amsfonts,latexsym}
\usepackage{color}
\date{ }

\begin{document}

\title{Regularity of Schr\"odinger's functional equation in the weak topology and moment measures
\thanks{2010 Mathematics Subject Classification : Primary 60G30 ; Secondary 93E20}
}

\author{Toshio Mikami\thanks{
%Department of Mathematics, Tsuda University,2-1-1 Tsuda-machi, Kodaira, Tokyo 187-8577, Japan,email: t.mikami@tsuda.ac.jp.
This work was  supported by JSPS KAKENHI Grant Numbers JP26400136 and JP16H03948.}\\
}

\maketitle

\begin{abstract}
We study the continuity and the measurability of the solution to Schr\"odinger's functional equation, with respect to space, kernel and marginals,
provided the space of all Borel probability measures is endowed with the weak topology.
This is a continuation of our previous result where the space of all Borel probability measures was endowed with the strong topology.
As an application, we construct a convex function of which the moment measure is  a given probability measure,
by the zero noise limit of a class of stochastic optimal transportation problems.
\end{abstract}

\maketitle

\section{Introduction}
E. Schr\"odinger considered the following problem to find the statistical property of a particle on a finite time interval. 
Suppose that there exist $N\ge 2$ particles in a set $A:=\{a_1, \cdots, a_{n_0}\}\subset {\bf R}^3$ and each particle moves independently, with a given transition probability, to a set $B:=\{b_1, \cdots, b_{n_1}\}\subset {\bf R}^3$, 
where   $1\le n_0, n_1\le N$.
Find the maximal probability of such events, provided the numbers of particles in each point $a_i$, $b_j$ are fixed (see section 7 in \cite{S2}  and also \cite{S1}).
Though he did not succeed in finding the maximal probability, he obtained Euler's equation for the variational problem above.
The continuum limit is called Schr\"odinger's functional equation (see \cite{B, C, Fo, J1} for the solution of this problem).
S. Bernstein \cite{S. B} generalized Schr\"odinger's idea and introduced the so-called Bernstein processes which are also called reciprocal processes.
The theory of stochastic differential equation for Schr\"odinger's functional equation was given by B. Jamison \cite{J2}.
The solution is Doob's h-path process  (see \cite{Doob}) with given two end point marginals.
Schr\"odinger's problem is also related to the theory of large deviations, the optimal mass transportation problem,
entropic estimates and functional inequalities
 (see, e.g. \cite{ADPZ, BCGL, Con0, Con, DLR, EMR, GLRT, Leo1, Leo2, 19,  21-1, 22, PAL, RR, R, V} and the references therein).
 
We describe E. Schr\"odinger's functional equation (see e.g. \cite{J1}) in the setting considered in this paper.
Let $S$ be a $\sigma$-compact metric space and $q\in C(S\times S; (0,\infty ))$.
%Schr\"odinger's functional equation can be described as follows.
For Borel probability measures $\mu_1, \mu_2$ on $S$,  find nonnegative $\sigma$-finite Borel measures $\nu_1, \nu_2$ on $S$ for which
the following holds:
\begin{equation}\label{11}
\begin{cases}
\displaystyle\mu_1(dx)=\nu_1(dx)\int_{S} q(x,y)\nu_2(dy),&\\
&\\
\displaystyle\mu_2(dy)=\nu_2(dy)\int_{S} q(x,y)\nu_1(dx).&
\end{cases} 
\end{equation} 
It is known that there exists a solution $(\nu_1,\nu_2)$ of (\ref{11})  (see \cite{C, J1}).
$(\nu_1,\nu_2)$ is unique up to a constant though the product measure $\nu_1\times\nu_2$ is unique.
Indeed,   for any $C>0$, $(C\nu_1,C^{-1}\nu_2)$ is also a solution of  (\ref{11}).
By the uniqueness of the solution to (1.1), we mean that of the product measure $\nu_1\times\nu_2$.
Let $\{K_m\}_{m\ge 1}$ be a nondecreasing sequence of compact subsets of $S$ such that $S=\cup_{m\ge 1} K_m$, where $K_1\equiv S$ when $S$ is compact.
When we consider $\nu_1$ and $\nu_2$ separately, considering $(C(\mu_1,\mu_2)\nu_1,C(\mu_1,\mu_2)^{-1}\nu_2)$ if necessary,
we assume that  the following holds:
\begin{equation}\label{17}
\nu_1 (K_{m (\mu_1,\mu_2)})=\nu_2(K_{m(\mu_1,\mu_2)}),
\end{equation}
where
$$m(\mu_1,\mu_2):=\min \{m\ge 1| \mu_1(K_m)\mu_2(K_m)>0\},  \quad 
C(\mu_1,\mu_2):=\left(\frac{\nu_2 (K_{m (\mu_1,\mu_2)})}{\nu_1(K_{m(\mu_1,\mu_2)})}\right)^{1/2}.$$
\begin{equation}\label{12}
\mu(dxdy):=\nu_1(dx)q(x,y)\nu_2(dy),
\end{equation} 
\begin{equation}\label{13}
u_i(x_i):=\log\biggl(\int_{S} q(x_1,x_2)\nu_j(dx_j)\biggr),\quad i,j =1,2,i\ne j.
\end{equation}
Then $\exp (u_1(x))$ and $\exp (u_2 (x))$ are positive %,$\nu_i(dx)=\exp (-u_i (x))\mu_i (dx)$ 
and 
\begin{equation}\label{14}
\mu(dxdy)=q(x,y)\exp (-u_1 (x)-u_2 (y))\mu_1 (dx) \mu_2(dy).
\end{equation}
(\ref{11}) can be rewritten as follows: for $i,j=1,2$, $i\ne j$,
\begin{equation}\label{15}
\mu_i (dx_i)=\exp (-u_i (x_i))\mu_i (dx_i)\int_{S} q(x_1,x_2)\exp (-u_j(x_j))\mu_j(dx_j).
\end{equation}
In particular, Schr\"odinger's problem (\ref{11}) is equivalent to finding functions $u_1$ and $u_2$
for which (\ref{15}) holds.
$(u_1,u_2)$ is unique up to a constant though $u_1 (x_1)+u_2(x_2)$ is unique.

Let ${\cal M}(S)$ and ${\cal P}(S)$ denote the space of all Radon measures and that of all Borel probability measures on $S$, respectively, 
where  a Radon measure means a locally finite and inner regular Borel measure.
It is easy to see that $\nu_1$ and $\nu_2$ are functionals of $\mu_1$, $\mu_2$ and $q$:
\begin{equation}\label{16}
\nu_i (dx)=\nu_i (dx;q,\mu_1, \mu_2),\quad u_i (x)=u_i (x;q,\mu_1, \mu_2),\quad i=1,2.
\end{equation}
In \cite{24}, we  considered the case where ${\cal P}(S)$ is endowed with
the strong topology and showed that if $S$ is compact, then the following is continuous:
$$\nu_i(dx;\cdot,\cdot,\cdot):C(S\times S)\times {\cal P}(S)\times {\cal P}(S) 
\longrightarrow{\cal M}(S),$$
$$\{u_i (x;\cdot,\cdot,\cdot)\}_{x\in S}:C(S\times S)\times {\cal P}(S)\times {\cal P}(S)
\longrightarrow C(S)$$
and $u_i\in C(S\times C(S\times S)\times {\cal P}(S)\times {\cal P}(S))$.
Here ${\cal M}(S)$ is endowed with the strong topology and 
$C(S\times S)$ and $C(S)$ are endowed with the topology induced by the uniform convergence
on $S\times S$ and $S$, respectively.
We also showed that if $S$ is $\sigma$-compact, then
the following is Borel measurable:
$$\int_S f(x)\nu_i(dx;\cdot,\cdot,\cdot):C(S\times S)\times {\cal P}(S)\times {\cal P}(S)
\longrightarrow\mathbb{R},\quad f\in C_0 (S)$$
$$u_i:S\times C(S\times S)\times {\cal P}(S)\times {\cal P}(S)\longrightarrow\mathbb{R}
\cup\{\infty\}.$$
As an application of this measurability result,
we showed that the coefficients of the mean field PDE system for
the h-path process with given two end point marginals
are measurable functions of space, time and marginal.
\begin{rk}
(\ref{17}) was assumed in \cite{24} and implies that for $i, j=1,2$,
$$\nu_i(dx_i;q,\mu_1,\mu_2)=
\frac{\nu_i(dx_i;q,\mu_1,\mu_2)\nu_j (K_{m(\mu_1,\mu_2)};q,\mu_1,\mu_2)}{(\nu_1(K_{m(\mu_1,\mu_2)};q,\mu_1,\mu_2)\nu_2(K_{m(\mu_1,\mu_2)};q,\mu_1,\mu_2))^{1/2}}.
$$
In particular, the measurability of $(q,\mu_1,\mu_2)\mapsto\nu_1(dx_1;q,\mu_1,\mu_2)\nu_2(dx_2;q,\mu_1,\mu_2)$
implies that of $(q,\mu_1,\mu_2)\mapsto\nu_i(dx_i;q,\mu_1,\mu_2)$.
\end{rk}

In this paper we consider the case where ${\cal P}(S)$ is endowed with the weak topology
and show the continuity and measurability results on $\nu_i$ and $u_i$ (see Theorem 2.1 and Corollaries 2.1-2.3 in section 2).
Our continuity result in the weak topology is useful when one considers the existence of a minimizer of a variational problem (see \cite{CDS} for the continuity result on optimal transport).
Indeed, it is not easy to show that a minimizing sequence is compact in the strong topology.
As an application (see Theorem 2.2 in section 2), we give a stochastic optimal transportation approach to moment measures (see \cite{Cor, Sa}).
The definition of a moment measure of a convex function is the following.
\begin{df}
Given a convex function $u :\mathbb{R}^d\longrightarrow\mathbb{R}\cup \{\infty\}$, the following is called the moment measure of $u$:
\begin{equation}\label{115}
\mu (dx):= (D u)_\# (\exp (-u(x))dx).
%P_1(dx)=(p_0 (x)dx)D(-\log p_0 (x))^{-1},
\end{equation}
\end{df}
\begin{rk}
If $\mu$ is a moment measure of  a convex function $u :\mathbb{R}^d\longrightarrow\mathbb{R}\cup \{\infty\}$,
then $\exp (-u(x))dx\delta_{D u(x)}(dy)$ is the unique minimizer of the $2$-Wasserstein distance $W_2 (\exp (-u(x))dx, \mu (dx))$,
provided $\int_{{\mathbb R}^d} \exp (-u(x))dx=1$ and $W_2 (\exp (-u(x))dx, \mu (dx))$ is finite  (see \cite{Bre1, Bre2, V}).
Here $\delta_x(dy)$ denotes the delta measure on $\{x\}$ and for $\mu_1,\mu_2\in {\cal P}(\mathbb{R}^d)$,
\begin{align}\label{1151}
W_2 (\mu_1,\mu_2):=&\biggl(\inf \biggl\{\int_{\mathbb{R}^d\times \mathbb{R}^d} |x-y|^2m(dxdy)\biggl|
m\in {\cal P}(\mathbb{R}^d\times \mathbb{R}^d),\\
&\qquad  m(dx\times \mathbb{R}^d)=\mu_1(dx), 
m(\mathbb{R}^d\times dy)=\mu_2(dy)\biggr\}\biggl)^{1/2}.\notag
\end{align}
\end{rk}

We describe an application of our continuity result more precisely.
Let $\varepsilon >0$ and let
 $W(t)$ and $\gamma (t)=\gamma (t;\omega )$ denote a $d$-dimensional Brownian motion and 
 a progressively measurable $\mathbb{R}^d$-valued stochastic process on a filtered probability space, respectively.
Consider the following SDE in a weak sense (see e.g. \cite{F}):
\begin{equation}\label{111}
dX^{\varepsilon, \gamma}  (t)=%X^{\varepsilon, \gamma}  (0)+\int_0^t 
\gamma (t)dt+\sqrt{\varepsilon}dW(t).
\end{equation}
For $P_0,P_1\in {\cal P}(\mathbb{R}^d)$,
\begin{equation}\label{112}
V_\varepsilon (P_0,P_1):=\inf \biggl\{ E\biggl[\int_0^1 \frac{1}{2\varepsilon }|\gamma (t)|^2dt \biggr]\biggr|
PX^{\varepsilon, \gamma} (t)^{-1}=P_t ,t=0,1 \biggr\},
\end{equation}
where $V_\varepsilon (P_0,P_1):=\infty$ if the set over which the infimum is taken is empty
(see \cite{ADPZ, BCGL, DLR, EMR, Jin, Leo2} for related problems on large deviations).
For $P\in  {\cal P}(\mathbb{R}^d)$, 
\begin{equation}\label{113}
{\cal S}(P):=
\begin{cases}
\displaystyle\int_{\mathbb{R}^d} p(x)\log p(x) dx,&\hbox{ if }\displaystyle p(x):=\frac{P(dx)}{dx} \hbox{exists,}\\
\displaystyle\infty,&\hbox{otherwise.}
\end{cases} 
\end{equation}
For $\varepsilon,r>0$, $P_1\in {\cal P}(\mathbb{R}^d)$,
\begin{align}\label{114}
&\Psi_{\varepsilon,r} (P_1)\\
:=&
\inf\biggl\{{\cal S}(P)-\varepsilon V_\varepsilon (P,P_1)+\frac{1}{2} \int_{\mathbb{R}^d}|x|^2  P(dx)\biggr| P(dx)=p(x)dx\in  {\cal P}(B_r)\biggr\},\notag
\end{align}
where 
$$B_r:=\{ x\in \mathbb{R}^d: |x|\le r\}.$$
By our weak continuity result of $(q,\mu_1,\mu_2)\mapsto \mu (dxdy;q,\mu_1,\mu_2)$,
we can easily prove the existence of a minimizer $P_{0,r,\varepsilon}$ of $\Psi_{\varepsilon,r} (P_1)$
from the lower semicontinuities of a relative entropy and  of ${\cal S}$ with respect to the weak topology (see (\ref{119}) and also Lemmas 3.4 and 3.5 in section 3).
We show that a subsequence of $\{p_{0,r,\varepsilon} (x)dx\}_{\varepsilon >0}$ weakly converges, as $\varepsilon\to 0$, to a Borel probability measure $p_0 (x)dx$
such that $-\log p_0 (x)$ is convex and $P_1$ is a moment measure of  $-\log p_0 (x)$.
This is formally implied by the representation of $P_{0,r,\varepsilon}$ and the SDE for the minimizer of 
$V_\varepsilon (P_{0,r,\varepsilon},P_1)$ (see  (\ref{29}) and (\ref{19})).
%Our approach can not be applied for the limit of $\mu (dxdy;g_\varepsilon (1),P_{0,r,\varepsilon}, P_1)$
%as $\varepsilon\to 0$ since it is not absolutely continuous with respect to the product measure of its marginals
%(see Remark 1.1).
We also show that $p_{0,r,\varepsilon} (x)$ has  a subsequence which uniformly converges, as $\varepsilon\to 0$, to $p_0 (x)$, provided $P_1$ is compactly supported. 

$\Psi_{\varepsilon,r} (P_1)$ formally converges, as $\varepsilon\to 0$, to the functional considered in \cite{Sa}
where they take the infimum over ${\cal P}(\mathbb{R}^d)$ instead of ${\cal P}(B_r)$.
Our approach makes the proof easier than \cite{Sa} since $ {\cal P}(B_r)$ is compact in the weak topology 
but can not be applied if we replace ${\cal P}(B_r)$ by ${\cal P}(\mathbb{R}^d)$, which we regret.

In the proof of the representation of $P_{0,r,\varepsilon}$ in (\ref{29}), we also make use of properties of the solution to Schr\"odinger's functional equation
and the duality theorem for $V_\varepsilon (P_0,P_1)$:
\begin{align}\label{116}
V_\varepsilon (P_0,P_1)
=&\sup\biggl\{\int_{\mathbb{R}^d} f (x)P_1(dx)-\int_{\mathbb{R}^d}\varphi  (0,x;f)P_0(dx)\biggr|f\in C_b^\infty (\mathbb{R}^d)\biggr\}.
\end{align}
Here the supremum is taken over all classical solutions $\varphi (t,x;f)$ to the following Hamilton-Jacobi-Bellman 
PDE:
\begin{align}\label{117}
\frac{\partial\varphi (t,x)}{\partial t}+\frac{\varepsilon}{2}\triangle_x \varphi (t,x)+\frac{\varepsilon}{2}|D_x\varphi  (t,x)|^2&=0,\quad 
(t,x)\in [0,1)\times \mathbb{R}^d,\\
\varphi (1,x) &=f(x)\notag
\end{align}
(see \cite{21-1,22, 25, TT} and the references therein).
\begin{align}\label{18}
g_\varepsilon (t,z):=&\frac{1}{(2\pi \varepsilon t)^{d/2}}\exp \left(-\frac{|z|^2}{2\varepsilon t}\right),
\quad t>0, z\in \mathbb{R}^d, \\
g_\varepsilon (t)(x,y):=&g_\varepsilon (t,y-x),\quad t>0, x,y\in \mathbb{R}^d.\notag
\end{align}
It is known that for any $P_0, P_1\in {\cal P}({\bf R}^d)$ for which $P_1 (dy)\ll dy$, 
there exists the unique weak solution to the following two end points problem of SDE (see \cite{J2} and also \cite{24, 25}):
\begin{align}\label{19}
dX^\varepsilon(t)=&\varepsilon D_x u_1 (X^\varepsilon(t);g_\varepsilon (1-t),PX^\varepsilon(t)^{-1}, P_1)dt
+\sqrt{\varepsilon} dW(t),\quad 0< t< 1,\\
PX^\varepsilon (t)^{-1}=&P_t,\quad t=0,1.\notag
\end{align}
$X^\varepsilon(t)$ is called the h-path process for $\sqrt{\varepsilon}  W(t)$ on $[0,1]$
with initial and terminal distribution $P_0$ and $P_1$, respectively.
The following is also known:
\begin{align}\label{110}
P(X^\varepsilon (0), X^\varepsilon (1))^{-1}(dxdy)
=&\nu_1(dx;g_\varepsilon (1), P_0,P_1)g_\varepsilon (1,y-x)\nu_2(dy;g_\varepsilon (1), P_0,P_1).
\end{align}
Suppose that $V_\varepsilon (P_0,P_1)$ is finite (see Remark 2.2 in section 2 for a sufficient condition).
Then $X^\varepsilon$ in (\ref{19}) is the unique minimizer of $V_\varepsilon (P_0,P_1)$ (see \cite{7, 11}, \cite{Leo1}-\cite{28}, \cite{TT}, \cite{Z} and the references therein).
Besides, there exists $f_o\in L^1 (P_1)$ which is unique up to a constant such that the following holds (see \cite{21-1,22, 24, 25, TT} and the references therein and also (\ref{14})):
\begin{align}\label{118}
f_o(y)-\varphi  (0,x;f_o)
=&\log p_1(y)-u_2 (y;g_\varepsilon (1),P_0, P_1)-u_1 (x;g_\varepsilon (1),P_0, P_1).
\end{align}
In particular, the following holds:
\begin{align}\label{119}
V_\varepsilon (P_0,P_1)
=&\int_{\mathbb{R}^d} f_0 (x)P_1(dx)-\int_{\mathbb{R}^d}\varphi  (0,x;f_0)P_0(dx)\\
=&{\cal S}(P_1)-\int_{\mathbb{R}^d} u_2 (x;g_\varepsilon (1),P_0, P_1)P_1(dx)\notag\\
&\qquad-\int_{\mathbb{R}^d}u_1 (x;g_\varepsilon (1),P_0, P_1)P_0(dx)\notag\\
=&H(P(X^\varepsilon (0), X^\varepsilon (1))^{-1} (dxdy)|P_0 (dx)g_\varepsilon (1,y-x) dy)\notag\\
%=&{\cal S}(P_1)+H(P(X^\varepsilon (0), X^\varepsilon (1))^{-1} (dxdy)|P_0 (dx)P_1 (dy))\notag\\
%&\qquad +\int_{\mathbb{R}^d\times \mathbb{R}^d}\log g_\varepsilon (1,y-x)P_0 (dx)P_1 (dy)\notag\\
=&{\cal S}(P_1)-H(P_0 (dx)P_1(dy)|P(X^\varepsilon (0), X^\varepsilon (1))^{-1} (dxdy))\notag\\
&\qquad-\int_{\mathbb{R}^d\times \mathbb{R}^d}\log g_\varepsilon (1,y-x)P_0 (dx)P_1 (dy).\notag
\end{align}
Here $H$ denotes the relative entropy of two measures: 
for $m, n \in {\cal P}(S\times S)$,
\begin{equation}\label{1201}
H(m|n)=
\begin{cases}
\int_{S\times S}\log \frac{m(dxdy)}{n(dxdy)}m(dxdy),& \hbox{ if } m\ll n,\\
\infty,&\hbox{ otherwise.}
\end{cases}
\end{equation}
%These facts also play crucial roles in the proof of our result.
\begin{rk}
If $V_\varepsilon (P_0,P_1)$ is finite, then $P_1 (dy)\ll dy$.
Indeed, $V_\varepsilon (P_0,P_1)$ is the relative entropy of $P(X^\varepsilon)^{-1}$ with respect to 
$P_0\ast P(\sqrt \varepsilon W)^{-1}$ on $C([0,1]; \mathbb{R}^d)$ and 
$$P_0\ast P(\sqrt \varepsilon W(1))^{-1}(dy)=\left( \int_{\mathbb{R}^d}g_\varepsilon (1,y-x) P_0 (dx)\right)dy.
$$
Here $\ast $ denotes the convolution of two measures.
\end{rk}
In section 2 we state our main results and prove them in section 4 by lemmas which are given in section 3.

\section{Main result}

In this section we state our main results.
We first describe assumptions precisely.

\noindent
(A1) $S$ is a complete $\sigma$-compact metric space.

\noindent
(A1)' $S$ is a compact metric space.

\noindent
(A2) $q\in C(S\times S ;(0,\infty))$.

We remark that ${\cal P}(S)$ is endowed with the weak topology 
and $C(S\times S)$ is endowed with  the topology induced by the uniform convergence on every compact subset of $S$.

Under (A1), let %$\{K_m\}_{m\ge 1}$ and 
$\{\varphi_m \}_{m\ge 1}$ be %, respectively,  
a nondecreasing sequence %of compact subsets  of $S$ and that 
of functions in $C_0 (S;[0,1])$ such that  the following holds:
$$%S=\cup_{m\ge 1} K_m,\quad 
\varphi_m(x)=1, \quad x\in K_m, m\ge 1.$$
(see (\ref{17})).
If $S=\mathbb{R}^d$, then $K_m:=B_m$ and we assume that $\varphi_m\in C_0 (B_{m+1};[0,1])$.
For $i\ne j$, $i,j=1,2$,
%, where $$B_r:=\{ x\in \mathbb{R}^d: |x|\le r\},\quad r>0.$$
\begin{equation}\label{21}
u_{i |m} (x_i;q,\mu_{1},\mu_{2}):=\log\left(\int_{S} q (x_1,x_2)\varphi_m(x_j) \nu_{j}(dx_j;q,\mu_{1},\mu_{2})\right),
\end{equation}
provided the right hand side is well defined (see (\ref{16}) and also (\ref{13})).
\begin{equation}\label{211}
\mu(dxdy;q,\mu_1, \mu_2):=\nu_1 (dx;q,\mu_1, \mu_2)q(x,y)\nu_2 (dx;q,\mu_1, \mu_2).
\end{equation}
The following is the continuity result on $\nu_1\times \nu_2$, $\mu$ and $u_{i|m}$.

\begin{thm}
Suppose that (A1) and (A2) hold and that $q_n\in C(S \times S;(0,\infty))$, $\mu_i, \mu_{i,n}\in {\cal P}(S)$, 
$n\ge 1$,  $i=1,2$ and 
\begin{equation}\label{22}
\lim_{n\to\infty}q_n=q,\quad \hbox{locally uniformly},
\end{equation}
\begin{equation}\label{23}
\lim_{n\to\infty}\mu_{1,n}\times \mu_{2,n}=\mu_{1}\times \mu_{2},\quad \hbox{weakly.}
\end{equation}
%(see (2.1) for notation).
Then for any $f\in C_0 (S\times S)$,
\begin{align}\label{24}
&\lim_{n\to \infty} \int_{S\times S}f(x,y)\nu_1 (dx;q_n,\mu_{1,n},\mu_{2,n})\nu_2 (dy;q_n,\mu_{1,n},\mu_{2,n})\\
=&\int_{S\times S}f(x,y)\nu_1 (dx;q,\mu_{1},\mu_{2})\nu_2 (dy;q,\mu_{1},\mu_{2}).\notag
\end{align}
In particular, 
\begin{equation}\label{241}
\lim_{n\to\infty}\mu (dxdy;q_n,\mu_{1,n},\mu_{2,n})=\mu (dxdy;q,\mu_1,\mu_2),\quad \hbox{weakly.}
\end{equation}
For any $\{x_{i,n}\}_{n\ge 1}\subset S$ which converges, as $n\to\infty$, to $x_i\in S$, $i=1, 2$
and for sufficiently large $m\ge 1$,
\begin{align}\label{25}
&\lim_{n\to\infty}\sum_{i=1}^2 u_{i | m} (x_{i,n};q_n,\mu_{1,n},\mu_{2,n})=\sum_{i=1}^2 u_{i | m}(x_i;q,\mu_1, \mu_2).
\end{align}
\end{thm}

Since $(\mu_1,\mu_2)\mapsto m(\mu_1,\mu_2)$ is measurable,
Theorem 2.1 implies the following.

\begin{co}
Suppose that (A1) and (A2) hold.
Then the following are Borel measurable: for $i=1,2$,
$$\int_{S} f(x)\nu_i(dx;\cdot,\cdot, \cdot):C(S\times S)\times {\cal P}(S)\times {\cal P}(S)
\longrightarrow\mathbb{R},\quad f\in C_0 (S),$$
$$u_i:S\times C(S\times S)\times {\cal P}(S)\times {\cal P}(S)
\longrightarrow\mathbb{R}\cup\{\infty\}.$$
\end{co}

If $S$ is compact, then $\nu_1 (S)=\nu_2 (S)$ (see (\ref{17})).
This implies, from Theorem 2.1, the following of which the proof is omitted.
\begin{co}
Suppose that (A1)' and the assumption of Theorem 2.1 except (A1) hold.
Then the following holds: for $i=1,2$,
\begin{align*}
\lim_{n\to \infty}  \int_{S}f(x)\nu_i (dx ;q_n,\mu_{1,n},\mu_{2,n})=&\int_{S}f(x)\nu_i (dx ;q,\mu_{1},\mu_{2}),
\quad f\in C (S),%\quad\hbox{weakly},
\notag
\end{align*}
and for any $\{x_{n}\}_{n\ge 1}\subset S$ which converges, as $n\to\infty$, to $x\in S$, 
\begin{align*}
\lim_{n\to\infty}u_{i} (x_{n};q_n,\mu_{1,n},\mu_{2,n})&=u_{i}(x;q,\mu_1, \mu_2).
\end{align*}
\end{co}

A uniformly bounded sequence of  convex functions on a convex neighborhood $N_A$ of a convex subset $A$ of $\mathbb{R}^d$ is compact in $C(A)$, provided $dist (A, N_A^c)$ is positive (see e.g., \cite{Ba}, section 3.3).
We describe an additional assumption and state a stronger result than above, provided $S\subset  \mathbb{R}^d$.

%\noindent (A2.$r$) $q\in C(B_{r} \times B_r ;(0,\infty))$.

\noindent
(A3.$r$)
There exists $C_r>0$ for which $x\mapsto C_r|x|^2+\log q(x,y)$ and 
$y\mapsto C_r |y|^2+\log q(x,y)$ are convex on $B_r$ for any $y\in B_r$ and any $x\in B_r$, respectively.

\begin{rk}
If $\log q(x,y)$ has bounded second order partial derivatives on $B_{r}$, then (A3.$r$) holds.
\end{rk}

\begin{align}\label{26}
||f||_{\infty,r}:=\sup_{x\in B_r}|f(x)|, \quad f\in C(B_r).
%||f||_{\infty}:=\sup_{x\in \mathbb{R}^d}|f(x)|, \quad f\in C(\mathbb{R}^d),\notag
\end{align}
The following is a stronger convergence result than Corollary 2.2.

\begin{co}
Let $r>0$.
Suppose that (A3.$r$) 
and the assumptions of Corollary 2.2 with $S=B_r$ hold.
Then for any $r'<r$,
\begin{align}\label{27}
&\lim_{n\to\infty}\sum_{i=1}^2||u_{i} (\cdot;q_n,\mu_{1,n},\mu_{2,n})-u_{i} (\cdot;q,\mu_1, \mu_2)||_{\infty,r'} =0.
\end{align}
\end{co}

\begin{equation}\label{28}
 {\cal P}_p(\mathbb{R}^d):=\biggl\{P\in {\cal P}(\mathbb{R}^d)\biggr|\int_{\mathbb{R}^d}|x|^p P(dx)<\infty\biggr\},
 \quad p\ge 1.
\end{equation}
As an application of our regularity result, we show that
there exists a convex function of which the moment measure is a given probability measure.

\begin{thm}
For  any $P_1(dx)=p_1(x)dx\in {\cal P}_2(\mathbb{R}^d)$ for which ${\cal S}(P_1)$ is finite, 
there exists a minimizer of $\Psi_{\varepsilon,r} (P_1)$.
For any minimizer $P_{0,r,\varepsilon}(dx)=p_{0,r,\varepsilon}(x)dx$ of $\Psi_{\varepsilon,r} (P_1)$,
\begin{equation}\label{29}
p_{0,r,\varepsilon}(x)= \frac{1}{C_\varepsilon}I_{B_r}(x)\exp \biggr(
-\varepsilon u_1(x;g_\varepsilon (1) ,P_{0,r,\varepsilon},P_1)-\frac{1}{2} |x|^2\biggr ),
\end{equation}
where $C_\varepsilon$ is a normalizing constant.
Besides, there exists a subsequence of $p_{0,r,\varepsilon}(x)dx$ which weakly converges, as $\varepsilon\to 0$, to a probability measure $p_0(x)dx$ such that $p_1(x)dx$ is a moment measure of  $-\log p_0$.
Suppose, in addition, that $P_1$ is compactly supported.
Then there exists a subsequence of $p_{0,r,\varepsilon}(x)$ which uniformly converges, as $\varepsilon\to 0$, to a probability density function $p_0(x)$ such that $p_1(x)dx$ is a moment measure of  $-\log p_0$.
\end{thm}

\begin{rk}
If $P_0, P_1(dx)=p_1(x)dx\in {\cal P}_2(\mathbb{R}^d)$ and ${\cal S}(P_1)$ is finite, 
then $V_\varepsilon (P_0, P_1)$ is finite.
Indeed,  from the last equality of (\ref{119}),
$$V_\varepsilon (P_0,P_1)
\le {\cal S}(P_1)-\int_{\mathbb{R}^d\times \mathbb{R}^d} \log g_\varepsilon (1, y-x)P_0(dx) p_1(y)dy
$$
since, the relative entropy is nonnegative.
%from (\ref{15}), by Jensen's inequality,
%$$u_2 (y;g_\varepsilon (1),P_0, P_1)\ge\int_{\mathbb{R}^d} (\log g_\varepsilon (1, y-x)-u_1 (x;g_\varepsilon (1),P_0, P_1))P_0(dx).$$
\end{rk}

%\newpage

\section{Lemmas}
 
 In this section we state and prove lemmas.
When it is not confusing, we omit the dependence of $u_i, \nu_i$ on $q,\nu_1,\nu_2$.

\subsection{Lemmas for the proof of Theorem 2.1 and Corollary 2.3}
 
The following lemma will be used in the proof of Theorem 2.1.

\begin{lm}
Suppose that (A1) and (A2) hold.
Then for any $\mu_1, \mu_2\in {\cal P}(S)$, $\mu$ defined by (\ref{12}),% and sufficiently large $m\ge 1$,
\begin{align}\label{33012}
\min_{x,y\in K_m} q (x,y)^{-1}\mu (K_m\times K_m)
\le &\int_{S}\varphi_m(x)\nu_{1}(dx)\int_{S}\varphi_m(y)\nu_{2}(dy)\\
\le &\max_{x,y\in supp(\varphi_m)} q (x,y)^{-1}.\notag
\end{align}
%$m\mapsto u_{1|m} (x_1)+u_{2|m} (x_2)$ is nondecreasing and the following holds:
%\begin{align}\label{3301}
%&\min_{x,y\in supp(\varphi_m) }\frac{q(x_{1},y)q(x,x_{2})}{q (x,y)}
%\times \int_{S\times S}\varphi_m(x)\varphi_m(y)\mu (dxdy)\\
%\le&\exp (u_{1|m} (x_1)+u_{2|m} (x_2))\notag\\
%\le &\max_{x,y\in supp (\varphi_{m})}\frac{q(x_{1},y)q(x,x_{2})}{q (x,y)},
%\quad x_1, x_2\in S.\notag
%\end{align}
\end{lm}
(Proof) The proof is done by the following (see (\ref{12})): % and (\ref{21})):
$$\nu_{1}(dx)\nu_{2}(dy)=q (x,y)^{-1}\mu(dxdy).\Box
$$
%\begin{equation}\label{3302}
%\exp (u_{1|m} (x_1)+u_{2|m} (x_2))=\int_{S\times S}\frac{q(x_{1},y)q(x,x_{2})}{q (x,y)}\varphi_m(x)\varphi_m(y) \mu (dxdy),
%\end{equation}
%provided the right hand side is positive.
%$\Box$

For $r>0$ and $q\in C(B_r\times B_r;(0,\infty))$,
 \begin{align}\label{31}
m_{q,r}:=&\min\{q(x,y)|\hbox{ }|x|, |y|\le r\},\\
M_{q,r}:=&\max\{q(x,y)|\hbox{ }|x|, |y|\le r\}.\notag
\end{align}

Lemmas 3.2 and 3.3 will be used  to prove Corollary 2.3.
 
\begin{lm}{\rm (\cite{B}, p. 194)}
Let $r>0$.
Suppose that (A2) with $S=B_r$  holds.
Then, for any $\mu_1, \mu_2\in {\cal P}(B_r)$,
%there exists a unique pair of nonnegative finite measures $\nu_{1}$, $\nu_{2}$ on $B_r$ for which (\ref{11})
%with $S=B_r$ and 
the following holds (see (\ref{13}) for notation):
%\begin{equation}\label{32}
%\frac{1}{\sqrt{M_{q,r}}}\le\nu_{1}(B_r)=\nu_{2}(B_r)\le \frac{1}{\sqrt{m_{q,r}}}.
%\end{equation}
\begin{equation}\label{3201}
\frac{m_{q,r}}{\sqrt{M_{q,r}}}\le  \exp( u_i (x))
\le \frac{M_{q,r}}{\sqrt{m_{q,r}}}, \quad x\in B_r, i=1,2.
\end{equation}
\end{lm}

%We state lemmas which will be used in the proof of Corollary 2.3.

%For the sake of completeness, we prove the following lemma which will be used in the proofs of Lemma 3.6 and Theorem 2.2.
By the method of proving the convexity of a log moment generating function,
we obtain the following.

\begin{lm}
Let $C$  and $\nu\in {\cal M}(\mathbb{R}^d)$ be a convex subset of $\mathbb{R}^d$ and a nonnegative Radon measure, respectively.
Suppose that $C\ni x\mapsto f(x,y) $ is convex, $\nu (dy)${\rm\text-a.e.}.
Then $C\ni x\mapsto \log \int_{\mathbb{R}^d }\exp (f (x,y))\nu(dy)$ is convex.
\end{lm}
(Proof)
For $x,y\in C$ and $\lambda \in (0,1)$, by H\"older's inequality, 
\begin{align}\label{411}
&\int_{\mathbb{R}^d}\exp (f(\lambda x+(1-\lambda ) y, x_2))\nu (dx_2)\\
\le &\int_{\mathbb{R}^d}\exp( \lambda f (x,x_2)+(1-\lambda )  f(y,x_2))  \nu(dx_2)\notag\\
\le &\left(\int_{\mathbb{R}^d} \exp (f (x,x_2))\nu(dx_2)\right)^{\lambda }
\left(\int_{\mathbb{R}^d} \exp (f(y,x_2)) \nu(dx_2)\right)^{1-\lambda }.\Box\notag
\end{align}

\subsection{Lemmas for the proof of Theorem 2.2}
In this subsection, we prove lemmas for the proof of Theorem 2.2.
Lemma 3.3 will be also used in the proof of Theorem 2.2.
$$B_{{\cal P}_p(\mathbb{R}^d ),r}:=\left\{P\in {\cal P}_p(\mathbb{R}^d )\biggl|\int_{\mathbb{R}^d}|x|^pP(dx)\le r\right\}.$$
The lower semicontinuity of a relative entropy and the continuity result in Theorem 2.1
imply the following.

\begin{lm}
Suppose that Theorem 2.1 holds.
Then for any $r, \varepsilon >0$,
the following is lower-semicontinuous on 
$B_{{\cal P}_2(\mathbb{R}^d ),r}\times B_{{\cal P}_2(\mathbb{R}^d ),r}$ (see  (\ref{13}), (\ref{16}) and (\ref{18}) for notation):
\begin{equation}\label{33}
\mu_1\times \mu_2
\mapsto 
-V_\varepsilon (\mu_{1}, \mu_{2})+{\cal S}(\mu_2)+\sum_{i=1}^2\int_{\mathbb{R}^d}
\frac{|x|^2}{2\varepsilon}\mu_i(dx).
\end{equation}
\end{lm}
(Proof)
From (\ref{119}), 
\begin{align}
&-V_\varepsilon (\mu_{1}, \mu_{2})+{\cal S}(\mu_2)+\sum_{i=1}^2\int_{\mathbb{R}^d}
\frac{|x|^2}{2\varepsilon}\mu_i(dx)\\
=&H(\mu_{1} (dx)\mu_2(dy)|\mu(dxdy;g_\varepsilon (1),  \mu_{1}, \mu_{2}))\notag\\
&\qquad +\frac{1}{\varepsilon}\langle \int_{\mathbb{R}^d}x\mu_{1} (dx), \int_{\mathbb{R}^d}y\mu_2 (dy)\rangle
-\log (2\pi\varepsilon)^{d/2}\notag
\end{align}
(see (\ref{110}) and  (\ref{211}) for notation).
Since $(m,n)\mapsto H(m(dxdy)|n(dxdy))$ is lower semicontinuous (see \cite{DE}, Lemma 1.4.3), the proof is over from Theorem 2.1.
$\Box$

The following lemma can be proved by the lower semicontinuity of a relative entropy. 

\begin{lm}
For any $r>0$, ${\cal S}$ is  lower-semicontinuous on $B_{{\cal P}_2(\mathbb{R}^d ),r}$ in the weak topology.
\end{lm}
(Proof)
$$q(x):=\frac{(1+|x|)^{-d-1}}{\int_{\mathbb{R}^d}(1+|y|)^{-d-1}dy},\quad x\in \mathbb{R}^d.$$
The proof is done by the following: 
\begin{equation}
{\cal S}(P)=H(P(dx)|q(x)dx)+\int_{\mathbb{R}^d} \log q(x)P(dx)
\end{equation}
(see e.g.  \cite{DE}, Lemma 1.4.3).$\Box$

\begin{lm}
Suppose that (A1) and (A2) hold.
Then for any $\mu_1, \mu_2\in {\cal P}(S)$, $\mu$ defined by (\ref{12}) and sufficiently large $m\ge 1$,
$m\mapsto u_{i|m}$ is nondecreasing, $i=1,2$ and the following holds:
\begin{align}\label{3301}
&\min_{x,y\in supp(\varphi_m) }\frac{q(x_{1},y)q(x,x_{2})}{q (x,y)}
\times \int_{S\times S}\varphi_m(x)\varphi_m(y)\mu (dxdy)\\
\le&\exp (u_{1|m} (x_1)+u_{2|m} (x_2))\notag\\
\le &\max_{x,y\in supp (\varphi_{m})}\frac{q(x_{1},y)q(x,x_{2})}{q (x,y)},
\quad x_1, x_2\in S.\notag
\end{align}
\end{lm}
(Proof) The proof is done by the following (see (\ref{12}) and (\ref{21})):
\begin{equation}\label{3302}
\exp (u_{1|m} (x_1)+u_{2|m} (x_2))=\int_{S\times S}\frac{q(x_{1},y)q(x,x_{2})}{q (x,y)}\varphi_m(x)\varphi_m(y) \mu (dxdy),
\end{equation}
provided the right hand side is positive.
$\Box$

For $i=1,2, m\ge 1, \varepsilon >0, x\in \mathbb{R}^d,$
\begin{align}\label{326}
\overline u_{i,\varepsilon} (x):=&\varepsilon u_i (x;g_\varepsilon (1),P_{0, r,\varepsilon}, P_1)+\frac{1}{2}|x|^2,\\
\overline u_{i|m,\varepsilon } (x):=&\varepsilon u_{i | m} (x;g_\varepsilon (1),P_{0, r,\varepsilon}, P_1)+\frac{1}{2}|x|^2.\notag
\end{align}
In the following lemma, the boundedness of the set $B_r$ plays a crucial role.

\begin{lm}
For any $\varepsilon, r>0$ and $P_1(dx)=p_1(x)dx\in {\cal P}(\mathbb{R}^d)$,
\begin{equation}\label{319}
\Psi_{\varepsilon,r} (P_1)\le -\log \{{\rm Vol} (B_r)\}+\frac{1}{2} \int_{B_r}|x|^2\frac{dx}{{\rm Vol}(B_r)}.
\end{equation}
Suppose %that $P_1\in {\cal P}_2(\mathbb{R}^d)$ and ${\cal S}(P_1)$ is finite and 
that $P_{0, r,\varepsilon}$ in (\ref{29}) is a minimizer of $\Psi_{\varepsilon,r} (P_1)$.
Then for $y_0:=\int_{\mathbb{R}^d} xP_1(dx)$,
\begin{equation}\label{320}
\exp (-\varepsilon {\cal S}(P_1)-\int_{\mathbb{R}^d}\frac{1}{2}|x|^2P_1(dx)-\Psi_{\varepsilon,r} (P_1))
\le C_\varepsilon\exp (-\overline u_{2,\varepsilon} (y_0)).
\end{equation}
In particular, for any sequence $\{\varepsilon_n\}_{n\ge 1}$  which converges to $0$ as $n\to\infty$,
the set $\{x\in B_r| \liminf_{n\to\infty }(\overline u_{1,\varepsilon_n} (x)+\overline u_{2,\varepsilon_n} (y_0))<\infty\}$ 
has a positive Lebesgue measure,
provided $P_1\in {\cal P}_2(\mathbb{R}^d)$ and ${\cal S}(P_1)$ is finite.
\end{lm}
(Proof)
Let $p_{uni, r}$ denote the probability density function of the uniform distribution on $B_r$.
Then the following implies (\ref{319}):
\begin{equation}\label{321}
\Psi_{\varepsilon,r} (P_1)
\le {\cal S}(p_{uni, r}(x)dx)+\frac{1}{2} \int_{B_r}|x|^2\frac{dx}{{\rm Vol}(B_r)}.
\end{equation}
We prove  (\ref{320}).
We only have to consider the case where ${\cal S}(P_1)$ is finite and
$P_1\in {\cal P}_2(\mathbb{R}^d)$.
From (\ref{119}) and (\ref{29}), by Jensen's inequality, 
\begin{align}\label{322}
\Psi_{\varepsilon,r} (P_1)
=& {\cal S}(P_{0, r,\varepsilon}) -\varepsilon 
\biggl({\cal S}(P_1)-\int_{\mathbb{R}^d} u_2 (x;g_\varepsilon (1),P_{0, r,\varepsilon}, P_1)P_1(dx)\\
&\qquad-\int_{\mathbb{R}^d}u_1 (x;g_\varepsilon (1),P_{0, r,\varepsilon}, P_1)P_{0, r,\varepsilon}(dx)\biggr)
+\int_{\mathbb{R}^d}\frac{1}{2}|x|^2P_{0, r,\varepsilon}(dx)\notag\\
=&-\log C_\varepsilon-\varepsilon {\cal S}(P_1) 
+\int_{\mathbb{R}^d} \overline u_{2,\varepsilon} (x)P_1(dx)
-\int_{\mathbb{R}^d}\frac{1}{2}|x|^2P_1(dx)\notag\\
\ge&-\log C_\varepsilon -\varepsilon {\cal S}(P_1)+\overline u_{2,\varepsilon} (y_0)
-\int_{\mathbb{R}^d}\frac{1}{2}|x|^2P_1(dx).\notag
\end{align}
Indeed, one can show that $\overline u_{2,\varepsilon}$ is convex from Lemma 3.3
and that 
$\overline u_{2,\varepsilon}$ is finite and continuous on $\mathbb{R}^d$ since
$\nu_1 (dx;g_\varepsilon (1),P_{0, r,\varepsilon}, P_1)$ is a finite measure on $B_r$.
The last part of this lemma can be shown by Fatou's lemma from (\ref{3301}) in Lemma 3.6 and from the following:
for $m>r$,
$$\overline u_{1,\varepsilon} (x)+\overline u_{2,\varepsilon} (y_0)
\ge \overline u_{1|m,\varepsilon} (x)+\overline u_{2,\varepsilon} (y_0), \quad \overline u_{2,\varepsilon}=\overline u_{2|m,\varepsilon}$$
since $\nu_1 (dx;g_\varepsilon (1),P_{0, r,\varepsilon}, P_1)$ is supported on $B_r$.$\Box$

For a convex function $f:\mathbb{R}^d\longrightarrow\mathbb{R}\cup\{\infty\}$, the $0$-sublevel set $f^{-1}((-\infty,0])$ is convex.
Roughy speaking, the following lemma can be proved from the fact that a uniformly bounded sequence of convex functions defined on the same open set is compact in the sup norm on any compact subset of the open set
(see section 3.3 in \cite{Ba}).

\begin{lm}
(i) For a convex set $C\subset \mathbb{R}^d$,
$dist (x,C)$ is a convex function.
(ii) For a bounded sequence of convex sets $\{C_n\subset \mathbb{R}^d\}_{n\ge 1}$,
there exists a closed convex set $C_\infty$ and 
a subsequence $\{C_{n_k}\}_{k\ge 1}$ of $\{C_n\}_{n\ge 1}$ such that
$\{dist (x,C_{n_k})\}_{k\ge 1}$  converges, as $k\to\infty$, to $dist (x,C_\infty)$
uniformly on every compact subset of $\mathbb{R}^d$.
(iii) For any $\gamma >0$, the following holds: for sufficiently large $k\ge 1$,
$$U_{-\gamma } (C_\infty):=\{y\in C_\infty |U_\gamma (y)\subset C_\infty\}\subset C_{n_k},$$
where $U_\gamma (y):=\{x\in \mathbb{R}^d: |x-y|<\gamma\}$.
\end{lm}
(Proof)
(i) For $x_1,x_2\in\mathbb{R}^d$, $\lambda\in (0,1)$, $y_1,y_2\in C$,
since $\lambda y_1+(1-\lambda) y_2\in C$,
\begin{align}
&dist (\lambda x_1+(1-\lambda) x_2,C)\\
\le& |\lambda x_1+(1-\lambda) x_2-(\lambda y_1+(1-\lambda) y_2)|\notag\\
\le& \lambda |x_1-y_1|+(1-\lambda) |x_2-y_2|.\notag
\end{align}
Taking the infimum over all  $y_1,y_2\in C$, the proof is done.

\noindent
(ii) Since $\{C_n\}_{n\ge 1}$  is bounded, $\{dist (x,C_n)\}_{n\ge 1}$ is also locally bounded, which implies that 
there exists a convex function $h(x)$ and a subsequence $\{dist (x,C_{n_k})\}_{k\ge 1}$ 
such that 
$$dist (x,C_{n_k})\to h(x),\quad k\to\infty,$$
uniformly on every compact subset of $\mathbb{R}^d$ (see, e.g., \cite{Ba}, section 3.3).
$$C_\infty:=h^{-1}(0).$$
Then it is easy to see that the set  $C_\infty$ is a closed convex set and $h(x)=dist (x,C_\infty)$.

\noindent
(iii) We only have to consider the case where $U_{-\gamma } (C_\infty)\ne\emptyset$.
From (ii), for sufficiently large $k\ge 1$,
\begin{equation}\label{332}
C_\infty\subset U_{\gamma_k} (C_{n_k}),
\end{equation}
where
$$\gamma_{k}:=\sup\left\{| dist (x,C_{n_k})|+\frac{\gamma}{2}\biggr|x\in C_\infty\right\}\to \frac{\gamma}{2}, \quad k\to\infty.$$
For $x\in U_{-\gamma } (C_\infty)$, if $x\notin C_{n_k}$, then the following which contradicts  (\ref{332}) holds: for $\tilde\gamma<\gamma$, 
$$\emptyset\ne U_{\gamma } (x)\cap U_{\tilde\gamma} (C_{n_k})^c\subset C_\infty\cap U_{\tilde\gamma} (C_{n_k})^c.$$
Indeed, since $C_{n_k}$ is convex, for $x\notin C_{n_k}$, there exists $p\in\mathbb{R}^d$ such that
$$\{y| \langle p, y-x\rangle\ge 0\}\subset C_{n_k}^c.\Box$$

\section{Proof of main results}
 
 In this section we  prove our main results.
 
 \noindent
(Proof of Theorem 2.1)
We first prove (\ref{24}).
For the sake of simplicity,
\begin{align}\label{41}
\nu_{i,n}(dx):=&\nu_i (dx;q_n,\mu_{1,n},\mu_{2,n}),\\ 
\mu_n (dxdy):=&\nu_{1,n}(dx)q_n (x,y)\nu_{2,n}(dy).\notag
\end{align}
Since $\{\mu_{1,n}(dx )=\mu_n (dx\times S),\mu_{2,n}(dy )=\mu_n (S\times dy)\}_{n\ge 1}$
is convergent, $\{\mu_n\}_{n\ge 1}$ is tight.
Indeed, for any Borel sets $A, B\in \mathbb{R}^d$,
$$\mu_n((A\times B)^c)\le\mu_n (A^c\times S)+\mu_n (S\times B^c)=\mu_{1,n} (A^c)+\mu_{2,n} (B^c),$$
and a convergent sequence of probability measures on a complete separable metric space is tight by Prohorov's Theorem (see, e.g., \cite{Bi}).
Here notice that a $\sigma$-compact metric space is separable.
By Prohorov's theorem, take a weakly convergent subsequence $\{\mu_{n_k}\}_{k\ge 1}$ and denote the limit by $\mu$.
Then it is easy to see that the following holds:
$$\mu_{1}(dx )=\mu (dx \times S),\quad \mu_{2}(dy )=\mu (S\times dy).$$
From (A2) and (\ref{22})-(\ref{23}), the following holds: for any $f\in C_0 (S\times S)$,
\begin{equation}\label{42}
\lim_{k\to\infty}\int_{S\times S}f(x,y)\nu_{1,n_k}(dx)\nu_{2,n_k}(dy)=
\int_{S\times S}f(x,y) q(x,y)^{-1}\mu (dxdy).
\end{equation} 
Indeed,
$$\nu_{1,n}(dx)\nu_{2,n}(dy)
=(\frac{1}{q_n (x,y)}-\frac{1}{q (x,y)})\mu_n (dxdy)
+\frac{1}{q (x,y)}\mu_n (dxdy).$$
The rest of the proof of (\ref{24}) is divided into the following (\ref{43})-(\ref{44}) which will be proved later.

\noindent
There exists a subsequence $\{\overline n_k\}\subset\{n_k\}$ and 
finite measures $\overline \nu_{1,m}$, $\overline \nu_{2,m}\in {\cal M}(supp (\varphi_m))$
such that for sufficiently large $m\ge 1$ and any $f\in C_0 (S\times S)$,
\begin{align}\label{43}
&\lim_{k\to\infty}\int_{S\times S}f(x,y)\varphi_m(x)\varphi_m(y)\nu_{1,\overline n_k}(dx)\nu_{2,\overline n_k}(dy)\\
=&\int_{S\times S}f(x,y)\overline \nu_{1,m}(dx)\overline \nu_{2,m}(dy).\notag
\end{align} 
From (\ref{43}), for sufficiently large $m\ge 1$ and any Borel sets $A_1,A_2\subset S$,
\begin{align}\label{44}
\int_{A_1\times A_2}q(x,y)^{-1}\mu (dxdy)
=&\frac{\int_{A_1\times K_m}q(x,y)^{-1}\mu (dxdy)\int_{K_m\times A_2}q(x,y)^{-1}\mu (dxdy)}
{\overline \nu_{1,m}(K_m)\overline \nu_{2,m}(K_m)}.
\end{align} 
(\ref{44}) implies that  $q(x,y)^{-1}\mu (dxdy)$ is a product measure which satisfies (\ref{11}).
(\ref{42}) and the uniqueness of the solution to (\ref{11}) implies that (\ref{24}) is true. 

\noindent
We prove (\ref{43})-(\ref{44}) to compete the proof of (\ref{24}).
(\ref{43}) can be proved by the diagonal method, since $\{\mu_n\}_{n\ge 1}$ is tight
and since for sufficiently large $m\ge 1$,
\begin{align}\label{45}
&\varphi_m(x_1)\varphi_m(x_2)\nu_{1,n_k}(dx_1)\nu_{2,n_k}(dx_2)\\
=&\int_{S}\varphi_m(x)\nu_{1,n_k}(dx)\int_{S}\varphi_m(y)\nu_{2,n_k}(dy)
\frac{\varphi_m(x_1)\nu_{1,n_k}(dx_1)}{\int_{S}\varphi_m(x)\nu_{1,n_k}(dx)}
\frac{\varphi_m(x_2)\nu_{2,n_k}(dx_2)}{\int_{S}\varphi_m(x)\nu_{2,n_k}(dx)}\notag
\end{align}
has a convergent subsequence from (\ref{33012}) in Lemma 3.1 by Prohorov's Theorem and since any weak limit is a product measure.
We prove (\ref{44}).
From (\ref{42}) and (\ref{43}), for sufficiently large $\tilde m\ge 1$,
\begin{align}\label{46}
&\int_{(A_1\times A_2)\cap (K_{\tilde m}\times K_{\tilde m})}q(x,y)^{-1}\mu (dxdy)\\
=&\int_{(A_1\times A_2)\cap (K_{\tilde m}\times K_{\tilde m})}\varphi_{\tilde m}(x)\varphi_{\tilde m}(y)q(x,y)^{-1}\mu (dxdy)\notag\\
=&\int_{(A_1\times A_2)\cap (K_{\tilde m}\times K_{\tilde m})}\overline \nu_{1,\tilde m}(dx)\overline \nu_{2,\tilde m}(dy)\notag\\
=&\overline \nu_{1,\tilde m}(A_1\cap K_{\tilde m})\overline \nu_{2,\tilde m}(A_2\cap K_{\tilde m}).\notag
\end{align} 
From (\ref{46}), for $\tilde m\ge m$, setting $A_i=K_{m}$, 
\begin{align}\label{47}
\overline \nu_{1,\tilde m}(A_1\cap K_{\tilde m})
=&\frac{\int_{(A_1\times K_m)\cap (K_{\tilde m}\times K_{\tilde m})}q(x,y)^{-1}\mu (dxdy)}
{\overline \nu_{2,\tilde m}(K_m)},\\
\overline \nu_{2,\tilde m}(A_2\cap K_{\tilde m})
=&\frac{\int_{(K_m\times  A_2)\cap (K_{\tilde m}\times K_{\tilde m})}q(x,y)^{-1}\mu (dxdy)}
{\overline \nu_{1,\tilde m}(K_m)},\notag\\
\overline \nu_{1,\tilde m}(K_m)\overline \nu_{2,\tilde m}(K_m)
=&\overline \nu_{1,m}(K_m)\overline \nu_{2,m}(K_m)=\int_{K_{m}\times K_{m}}q(x,y)^{-1}\mu (dxdy).\notag
\end{align} 
Substitute (\ref{47}) to (\ref{46}) and let $\tilde m\to\infty$. Then we obtain  (\ref{44}).
(\ref{25}) can be shown from  (\ref{24}) by the following: from (\ref{21}),
\begin{align}
&\exp \left(\sum_{i=1}^2 u_{i|m} (x_{i,n};q_n,\mu_{1,n},\mu_{2,n})\right)\\
=&
\int_{S\times S}q_n(x_{1,n},y)q_n(x,x_{2,n})\varphi_m(x)\varphi_m(y) \nu_{1,n}(dx)\nu_{2,n}(dy),\notag
\end{align}
provided the right hand side is positive.$\Box$

For a compact set $K\subset S$, ${\cal P}(S)\ni\nu\mapsto \nu (K)$ is upper semicontinuous in the weak topology and is hence measurable.
Corollary 2.1 can be proved in the same way as in \cite{24} and we omit the proof.

As we mentioned in section 2, we omit the proof of Corollary 2.2.
Corollary 2.2 and Lemmas 3.2 and 3.3 immediately imply Corollary 2.3 (see \cite{Ba}, section 3.3)
and we omit the proof.
Indeed, if any subsequence of a sequence of pointwise convergent continuous functions  has a uniformly convergent subsequence,
then it is uniformly convergent.

Before we prove Theorem 2.2, we briefly describe the idea of the proof.
Theorem 2.1 and Lemmas 3.4 - 3.5 imply the lower semicontinuity of the functional that we minimize in $\Psi_{\varepsilon, r}(P_1)$.
(\ref{319}) in Lemma 3.7 implies the finiteness of $\Psi_{\varepsilon, r}(P_1)$.
In particular, the existence of a minimizer $p_{0,r,\varepsilon}$ of $\Psi_{\varepsilon, r}(P_1)$ is obtained.
(\ref{29}) can be proved by the Duality Theorem (\ref{116}) for $V_\varepsilon (P_{0,r,\varepsilon}, P_1)$
and by 
the fact that the relative entropy of two probability measures is nonnegative and is equal to zero if and only if two probability measures are the same.
The characterization of the limit $p_{0}$ of $p_{0,r,\varepsilon}$, as $\varepsilon\to 0$, 
can be inferred from the following.
Roughly speaking, from \cite{19},
$$(2\varepsilon V_\varepsilon (P_{0,r,\varepsilon}, P_1))^{1/2}\sim W_2 (p_{0}(x)dx, P_1), \quad \varepsilon\to 0$$ 
(see (\ref{1151}) and (\ref{112}) for notation).
Besides, there exists a convex function $u:\mathbb{R}^d\longrightarrow\mathbb{R}\cup\{\infty\}$ such that 
for the minimizer  $X^\varepsilon$ of $V_\varepsilon (P_{0,r,\varepsilon}, P_1)$,  
as $\varepsilon\to 0$,
$$\varepsilon D_x u_1 (x;g_\varepsilon (1-t),PX^\varepsilon (t)^{-1}, P_1)+x\sim Du(x),\quad 0\le t\le 1,$$
$$\mu(dxdy;g_\varepsilon (1),P_{0,r,\varepsilon}, P_1)\sim P_0(dx) \delta_{D u (x)}(dy).$$
In particular, 
$$\varepsilon   u_1 (x;g_\varepsilon (1),P_{0,r,\varepsilon}, P_1)+\frac{|x|^2}{2} \sim u(x)+
{\rm Constant}.$$

(Proof of Theorem 2.2)
Since ${\cal P}(B_r)$ is tight, by Prohorov's Theorem (see, e.g., \cite{Bi}),
Lemmas 3.4-3.5 and (\ref{319}) in Lemma 3.7
imply the existence of a minimizer 
$P_{0,r,\varepsilon}(dx)=p_{0,r,\varepsilon}(x) dx$ of $\Psi_{\varepsilon,r} (P_1)$.
By (\ref{116}),
\begin{align}\label{413}
\Psi_{\varepsilon,r} (P_1)=&
\inf\biggl\{{\cal S}(p(x) dx)-\varepsilon \left(\int_{\mathbb{R}^d}f(x)p_1(x)dx-\int_{\mathbb{R}^d}\varphi(0,x;f)p(x)dx\right)\\
&\qquad +\frac{1}{2} \int_{\mathbb{R}^d}|x|^2 p(x)dx
\biggr| p(x)dx\in  {\cal P}(B_r), f\in C_b^\infty (\mathbb{R}^d)\biggr\}.\notag
\end{align}
Let $f_{0,r,\varepsilon}$ denote $f_o$ in (\ref{118}) with $P_0=P_{0,r,\varepsilon}$.
Then\begin{align}\label{414}
\Psi_{\varepsilon,r} (P_1)=&
\inf\biggl\{{\cal S}( p(x) dx)+\int_{\mathbb{R}^d}\left(\varepsilon \varphi(0,x;f_{0,r,\varepsilon}) +\frac{|x|^2}{2} \right)p(x)dx\\
&\qquad  \biggr|  p(x)dx\in  {\cal P}(B_r)\biggr\}
-\varepsilon \int_{\mathbb{R}^d}f_{0,r,\varepsilon}(x)p_1(x)dx\notag
\end{align}
(see (\ref{118}), (\ref{14}) and Remark 2.2).
Indeed,  for $p(x)dx\in  {\cal P}(B_r)$,
$$\int_{\mathbb{R}^d}\varphi(0,x;f_{0,r,\varepsilon}) p(x)dx-\int_{\mathbb{R}^d}f_{0,r,\varepsilon}(x)p_1(x)dx
\ge -V_\varepsilon (p(x)dx,P_1),$$
since 
\begin{align*}
&\int_{\mathbb{R}^d}\varphi(0,x;f_{0,r,\varepsilon}) p(x)dx-\int_{\mathbb{R}^d}f_{0,r,\varepsilon}(x)p_1(x)dx\\
=&-\int_{\mathbb{R}^d\times \mathbb{R}^d}
\log\left (\frac{\mu(dxdy;g_\varepsilon (1),P_{0,r,\varepsilon}, P_1)}{P_{0,r,\varepsilon}(dx)g_\varepsilon (1,y-x)}\right)\mu(dxdy;g_\varepsilon (1),p(x)dx, P_1),
\end{align*}
\begin{align*}
&V_\varepsilon (p(x)dx,P_1)\\
=&-\int_{\mathbb{R}^d\times \mathbb{R}^d}
\log\left (\frac{p(x)dxg_\varepsilon (1,y-x)}{\mu(dxdy;g_\varepsilon (1),p(x)dx, P_1)}\right)\mu(dxdy;g_\varepsilon (1),p(x)dx, P_1),
\end{align*}
and by Jensen's inequality, 
\begin{align*}
\int_{\mathbb{R}^d\times \mathbb{R}^d}
\log\left (\frac{\mu(dxdy;g_\varepsilon (1),P_{0,r,\varepsilon}, P_1)}{P_{0,r,\varepsilon}(dx)g_\varepsilon (1,y-x)}\frac{p(x)dxg_\varepsilon (1,y-x)}{\mu(dxdy;g_\varepsilon (1),p(x)dx, P_1)}\right)&\\
\qquad\times \mu(dxdy;g_\varepsilon (1),p(x)dx, P_1)&\le 0.
\end{align*}
(\ref{29}) holds since 
$\varphi(0,x;f_{0,r,\varepsilon})-u_1(x;g_\varepsilon (1) ,P_{0,r,\varepsilon},P_1)$
is a constant $C$ (see (\ref{118})) and since, for $p(x)dx\in  {\cal P}(B_r)$,
\begin{align*}
&S(p(x) dx)+\int_{\mathbb{R}^d}\left(\varepsilon \varphi(0,x;f_{0,r,\varepsilon}) +\frac{|x|^2}{2} \right)p(x)dx\\
=&\int_{B_r} p(x)dx \log\frac{p(x)}{C_\varepsilon^{-1}I_{B_r}(x)\exp (-\varepsilon u_1(x;g_\varepsilon (1) ,P_{0,r,\varepsilon},P_1)-\frac{1}{2} |x|^2 )}-\log C_\varepsilon +C\\
\ge &-\log C_\varepsilon +C.
\end{align*}
Here 
$$C_\varepsilon:=\int_{B_r} \exp \left(-\varepsilon u_1(x;g_\varepsilon (1) ,P_{0,r,\varepsilon},P_1)-\frac{1}{2} |x|^2 \right)dx,$$
and the equality holds if and only if 
$$p(x)=\frac{1}{C_\varepsilon}I_{B_r}(x)\exp\left (-\varepsilon u_1(x;g_\varepsilon (1) ,P_{0,r,\varepsilon},P_1)-\frac{1}{2} |x|^2 \right).$$
We prove the second part of Theorem 2.2.
For $i=1,2$, $m\ge 1$ and $x\in \mathbb{R}^d$,
\begin{align}\label{415}
\nu_{i,\varepsilon } (dx):=&\nu_{i} (dx;g_\varepsilon (1),  P_{0,r,\varepsilon}, P_1), \\
\mu_\varepsilon (dxdy):=&\nu_{1,\varepsilon } (dx)g_\varepsilon (1,y-x)\nu_{2,\varepsilon } (dy),\notag\\
u_{i|m, \varepsilon}(x):=&u_{i|m} (x;g_\varepsilon (1),  P_{0,r,\varepsilon}, P_1)\notag
\end{align}
(see (\ref{21}) for notation).
Since ${\cal P}(B_r)$ is compact, $\{P_{0,r,\varepsilon}\}_{\varepsilon >0}$ and $\{\mu_\varepsilon\}_{\varepsilon >0}$ has a weakly convergent subsequence by Prohorov's theorem in the same way as in the proof of Theorem 2.1 (see \cite{Bi}).
Let $P_0$ and $\mu$ denote the weak limit along the same subsequence, as $\varepsilon \to 0$,  of $P_{0,r,\varepsilon}$ and $\mu_\varepsilon$, respectively.
For sufficiently large  $m\ge 1$, by the diagonal method,
$\overline u_{1|m, \varepsilon} (x)+\overline u_{2|m, \varepsilon} (y)$ has a subsequence which is uniformly convergent, as $\varepsilon\to 0$, on every compact subset of $\mathbb{R}^d\times \mathbb{R}^d$ (see (\ref{326}) for notation).
Indeed, for sufficiently large  $m\ge 1$ and small $\varepsilon>0$,
$\overline u_{i|m, \varepsilon}$, $i=1,2$ are convex from Lemma 3.3,
and $\overline u_{1|m, \varepsilon} (x)+\overline u_{2|m, \varepsilon} (y)$ is uniformly bounded
on every compact subset of $\mathbb{R}^d\times \mathbb{R}^d$, from (\ref{3301}) in Lemma 3.6:
\begin{align}\label{416}
&|\overline u_{1|m, \varepsilon} (x)+\overline u_{2|m, \varepsilon} (y)+\varepsilon\log (2\pi\varepsilon)^{d/2}|\\
\le &(m+1)^2+(m+1)(|x|+|y|)-\varepsilon \log \mu_\varepsilon (B_{m}\times B_{m}),
\quad x, y\in \mathbb{R}^d.\notag
\end{align}
Let $\{\varepsilon_n\}_{n\ge 1}$ denote a sequence which converges to $0$, as $n\to\infty$ and along which the above sequences are all convergent.
\begin{equation}\label{4170}
u_m (x,y):=\lim_{n\to\infty} (\overline u_{1|m,  \varepsilon_n} (x)+\overline u_{2|m,  \varepsilon_n} (y)),
\quad m\ge 1, x, y\in \mathbb{R}^d.
\end{equation}
There exists the limit
\begin{equation}\label{417}
u (x,y):=\lim_{m\to\infty}u_m (x,y),\quad x, y\in \mathbb{R}^d.
\end{equation}
Indeed,  $m\mapsto u_m$ is nondecreasing since
\begin{align}\label{41711}
&\exp (u_{1|m, \varepsilon} (x_1)+u_{2|m,\varepsilon} (x_2))\\
=&\int_{S\times S}
\frac{g_\varepsilon (1,x_{1}-y)g_\varepsilon (1,x-x_{2})}{g_\varepsilon (1,x-y)}\varphi_m(x)\varphi_m(y) \mu_\varepsilon (dxdy).\notag
\end{align}
From the last statement of Lemma 3.7, there exists  $x_0\in B_r$ such that $u(x_0, y_0)<\infty$, since
$$\overline u_{1,  \varepsilon_n}(x)+\overline u_{2,  \varepsilon_n}(y_0)
\ge \overline u_{1|m,  \varepsilon_n} (x)+\overline u_{2|m,  \varepsilon_n} (y_0).$$
To complete the proof of Theorem 2.2, we show that the following holds:
\begin{equation}\label{419}
\langle x,y\rangle=u (x,y),\quad \mu{\rm\text-a.s.},
\end{equation}
\begin{equation}\label{420}
x=D_y u (x_0,y),\quad y=D_x u (x,y_0),\quad \mu{\rm\text-a.s.},
\end{equation}
\begin{equation}\label{421}
P_0(dx)=\frac{I_D (x)}{C}\exp (-u (x,y_0))dx,
\end{equation}
where $D$ is a convex  subset of $B_r$ and $C$ is a normalizing constant.
Notice that $u (x,y)$ is convex and is differentiable a.e. on its domain.

\noindent
$\underline {\hbox{Proof of (\ref{419})}}$
The following implies that (\ref{419}) holds:
for sufficiently large $m\ge r$, 
\begin{equation}\label{422}
\langle x,y\rangle=u_m (x,y), \quad   (x,y)\in \mathbb{R}^d\times Int (supp (\varphi_m)), \mu{\rm\text-a.s.}.
\end{equation}
Indeed,  from (\ref{417}) and (\ref{422}), 
for  sufficiently large $m>r$,
$$\langle x,y\rangle=u_m (x,y)=u_{m' }(x,y)=u (x,y),\quad m'>m,
$$
$(x,y)\in \mathbb{R}^d\times Int (supp (\varphi_m))(\subset \mathbb{R}^d\times Int (supp (\varphi_{m'})))$, $\mu${\rm\text-a.s.}.
To prove (\ref{422}), we first prove that the following holds:  for sufficiently large $m\ge r$, 
\begin{equation}\label{424}
\langle x,y\rangle\le u_m (x,y), \quad  (x,y)\in  supp (P_0)\times (supp (P_1)\cap Int (supp (\varphi_m))).
\end{equation}
For $i\ne j$, $i,j=1,2$,\begin{equation}\label{4241}
\mu_{i|m,\varepsilon}(dx_i):=\nu_{i,\varepsilon} (dx_i)\exp (u_{i|m,\varepsilon} (x_i))
=\int_{\{x_j\in \mathbb{R}^d\}} \varphi_m(x_j)\mu_{\varepsilon}(dx_1dx_2).
\end{equation}
Then for $\delta >0$ and $(x,y)\in  supp (P_0)\times (supp (P_1)\cap Int (supp (\varphi_m)))$,
\begin{align}\label{425}
&\exp \left(\frac{\overline u_{1|m,\varepsilon} (x)+\overline u_{2|m,\varepsilon} (y)}{\varepsilon}\right)\\
=&\int_{\mathbb{R}^d\times \mathbb{R}^d}\frac{1}{(2\pi \varepsilon)^d}
\exp  \left(\frac{\langle x_1,y\rangle+\langle x,x_2\rangle-\overline u_{1, \varepsilon|m } (x_1)-\overline  u_{2, \varepsilon|m } (x_2)}{\varepsilon}\right)\notag\\
&\qquad\times \varphi_m(x_1)\varphi_m(x_2)\mu_{1|m,\varepsilon}(dx_1)\mu_{2|m,\varepsilon}(dx_2)\notag\\
\ge &\int_{U_\delta (x)\times U_\delta (y)}\frac{1}{(2\pi \varepsilon)^d}
\exp  \left(\frac{\langle x_1,y\rangle+\langle x,x_2\rangle -\overline u_{1, \varepsilon|m} (x_1)-\overline  u_{2, \varepsilon|m } (x_2)}{\varepsilon}\right)\notag\\
&\qquad\times\varphi_m(x_2) \mu_{1|m,\varepsilon}(dx_1)p_1(x_2)dx_2\notag
\end{align}
(see (\ref{41711})).
Indeed, for $m>r$, $\mu_{1|m,\varepsilon}(dx)$ is supported on $B_r$ since $P_{0,r, \varepsilon}\in {\cal P}(B_r)$
and 
$$\mu_{2|m,\varepsilon }(dy)=\left (\int_{B_r} g_\varepsilon (1,y-x)\varphi_m(x)\nu_{1,\varepsilon } (dx)\right)\nu_{2,\varepsilon } (dy) =p_1(y)dy.
$$
(\ref{425}) implies (\ref{424}) since
\begin{align*}
&\int_{U_\delta (y)}\varphi_m(x_2)p_1(x_2)dx_2>0,\notag\\
&\liminf_{n\to \infty}\mu_{1|m,\varepsilon_n}(U_\delta (x))\ge
\int_{U_\delta (x)\times \mathbb{R}^d}\varphi_m(x_2)\mu(dx_1dx_2)\\
\ge &P_0(U_\delta (x))-P_1(B_m^c)>0,\quad \hbox{for  sufficiently large $m$}.\notag
\end{align*}
Next we prove that the following holds:  for sufficiently large $m\ge 1$, 
\begin{equation}\label{427}
\langle x,y\rangle\ge u_m (x,y), \quad  \mu{\rm\text-a.s.}. 
\end{equation}
$$A_{m,\delta, k}:=\{(x,y)\in B_r\times U_k (0) |\langle x,y\rangle-u_m (x,y)<-\delta\},\quad \delta> 0, k\ge 1.$$
Then $A_{m,\delta,k}$ is open since $u_m$ is convex and finite (see (\ref{416})-(\ref{4170})) and is continuous.
The following implies that (\ref{427})  is true: from (\ref{4170}),  for sufficiently large $m\ge 1$, 
\begin{align}\label{428}
\mu(A_{m,\delta,k})\le&\liminf_{ n\to \infty }\mu_{\varepsilon_n}(A_{m,\delta,k}),\\
\mu_{\varepsilon_n}(A_{m,\delta,k})=&\int_{A_{m,\delta,k}}\frac{1}{
(2\pi \varepsilon_n)^{d/2}}\exp \left(\frac{\langle x,y\rangle-\overline u_{1|m, \varepsilon_n} (x)-\overline  u_{2|m, \varepsilon_n } (y)}{\varepsilon_n}\right)\notag\\
&\qquad\times  \mu_{1|m,\varepsilon_n}(dx)\mu_{2|m,\varepsilon_n}(dy)\notag\\
\to &0\quad n\to \infty .\notag
\end{align}

\noindent
$\underline {\hbox{Proof of (\ref{420})}}$
For $(x,y)\in supp (P_0)\times supp (P_1)$,
\begin{equation}\label{429}
\langle x,y\rangle\le u (x,y)=u (x,y_0)+u (x_0,y)-u (x_0,y_0).
\end{equation}
Indeed, from  (\ref{417}) and (\ref{424}), 
for sufficiently large $m> r$ such that $y\in Int (supp (\varphi_m))$,
$$\langle x,y\rangle\le u_m (x,y)\le u (x,y).$$
(\ref{417}) and the following imply (\ref{429}): from  (\ref{4170}), 
$$u_m (x,y)=u_m (x,y_0)+u_m(x_0,y)-u_m (x_0,y_0).$$
$$A:=\{(x,y)\in \mathbb{R}^d\times\mathbb{R}^d |\langle x,y\rangle=u (x,y)\}.$$
$u (x,y_0)$ and $u (x_0,y)$ are finite for $(x,y)\in A$ from (\ref{416}), since
from (\ref{417}) and the equality in (\ref{429}),
\begin{align*}
&\langle x,y\rangle= u (x,y)\\
\ge &\max (u (x,y_0)+u_m (x_0,y)-u (x_0,y_0), u_m (x,y_0)+u (x_0,y)-u (x_0,y_0)).
\end{align*}
For a set $B\subset \mathbb{R}^d$ and a function $f:B\longrightarrow\mathbb{R}$,
$$
co\hbox{ }B:=\left\{\sum_{i=1}^{d+1}\lambda_i x_i\biggl | \sum_{i=1}^{d+1}\lambda_i=1, 
\lambda_i\ge 0, x_i\in B,1\le i\le d+1 \right\},
$$
$$con\hbox{ }f(x):=
\begin{cases}
\displaystyle\inf\biggl\{\sum_{i=1}^{d+1}\lambda_i f(x_i)\biggl | x=\sum_{i=1}^{d+1}\lambda_i x_i, 
\sum_{i=1}^{d+1}\lambda_i=1, &\\
\qquad \lambda_i\ge 0, x_i\in B,1\le i\le d+1\biggr\}, \quad  x\in co\hbox{ }B,&\\
\displaystyle\infty, \quad x\notin co\hbox{ }B.&
\end{cases} 
$$
Then, from (\ref{429}),  for $x\in  supp (P_0)$,
\begin{align}\label{430}
u (x,y_0)-u (x_0,y_0)\ge &\sup\{ \langle x,y\rangle-u (x_0,y)|y\in supp (P_1)\}\\
=&\sup\{ \langle x,y\rangle-con\hbox{ }(u|_{supp (P_1)}) (x_0,y)|y\in \mathbb{R}^d\}.\notag
\end{align}
Here $(u|_{supp (P_1)}) (x_0,y)$ denotes the restriction of $u (x_0,y)$ on $supp (P_1)$
and the equality holds if $(x,y_x)\in A$ for some $y_x\in supp (P_1)$, in which case $x\in \partial_y con\hbox{ }(u|_{supp (P_1)}) (x_0,y_x)$, where for a function $f:\mathbb{R}^d
\longrightarrow\mathbb{R}\cup\{\infty\}$,
$$\partial_y f(y)
:=\{x\in \mathbb{R}^d| f(z)\ge f(y)+\langle x,z-y\rangle, \hbox{ for any }z\in \mathbb{R}^d\}.
$$
In particular, $x\in \partial_y con\hbox{ }(u|_{supp (P_1)}) (x_0,y)$, $\mu${\rm\text-a.s.} from (\ref{419}).
$x=D_y u (x_0,y), \mu${\rm\text-a.s.}
since 
$$con\hbox{ }(u|_{supp (P_1)}) (x_0,y)=u (x_0,y),\quad y\in supp (P_1),$$
$$\partial_y con\hbox{ }(u|_{supp (P_1)}) (x_0,y)=\{D_y u (x_0,y)\},\quad dy{\rm\text-a.e.}
\hbox{ on $supp (P_1)$}$$
and since $P_1(dx)$ has a probability density function.
In the same way, one can show that $y=D_x u (x,y_0)$, $\mu${\rm\text-a.s.}.

\noindent
$\underline {\hbox{Proof of (\ref{421})}}$
$$D_{R,  \varepsilon }:=\{x\in B_r| \overline u_{1,\varepsilon}(x)+\overline u_{2,\varepsilon}(y_0)\le R\}, \quad R, \varepsilon >0$$
(see (\ref{326}) for notation).
Then, from Lemma 3.7, 
\begin{equation}\label{431}
\lim_{R\to \infty}\limsup_{\varepsilon\to 0}\int_{D_{R,  \varepsilon }^c}   p_{0,r,\varepsilon}(x)dx=0.
\end{equation}
Indeed,
\begin{align*}
\int_{D_{R,  \varepsilon }^c}   p_{0,r,\varepsilon}(x)dx
= &\frac{1}{C_\varepsilon \exp (-\overline u_{2,\varepsilon}(y_0))}
 \int_{B_r\cap D_{R,  \varepsilon }^c} \exp (-\overline u_{1,\varepsilon}(x)-\overline u_{2,\varepsilon}(y_0))dx\\
 \le & \exp (-R)\frac{{\rm Vol} (B_r)}{C_\varepsilon \exp (-\overline u_{2,\varepsilon}(y_0))}.\notag
\end{align*}
For $\delta >0$,  
$$\psi_{\delta, R, \varepsilon }(x):= \max \left(0, 1-\frac{dist (x,D_{R,  \varepsilon })}{\delta} \right).$$
Then, from Lemma 3.8, there exists a convergent subsequence $\{\psi_{\delta, R, \varepsilon_{n_k} }(x)\}_{k\ge 1}$ in $C(B_r)$
and a closed convex set $D_{R,0}\subset B_r$ such that 
\begin{equation}\label{4311}
\lim_{k\to\infty } ||\psi_{\delta, R, \varepsilon_{n_k} }-\psi_{\delta, R, 0}||_{\infty,r}=0.
\end{equation}
$$D:=\cup_{R>0}  D_{R, 0}.$$
Then we prove that the following holds: for a closed set $B\subset B_r$,
\begin{align}\label{4312}
&\lim_{R\to \infty}\limsup_{k\to\infty }\int_{B\cap D_{R,  \varepsilon_{n_k} }}  p_{0,r,\varepsilon_{n_k}}(x)dx\\
&\le \frac{1}{\int_{D} \exp (-u(x, y_0))dx}\int_{B\cap D} \exp (-u(x, y_0))dx.\notag
\end{align}
The proof of (\ref{4312}) is done by the following (\ref{4320})-(\ref{43201}) which will be proved later.
\begin{align}\label{4320}
&\lim_{R\to\infty}\limsup_{k\to\infty }\int_{B\cap D_{R,  \varepsilon_{n_k} }} \exp (-\overline u_{1,\varepsilon_{n_k}}(x)-\overline u_{2,\varepsilon_{n_k}}(y_0))dx\\
\le & \int_{B\cap D}\exp (-u(x,y_0))dx,\notag
\end{align}
\begin{equation}\label{43201}
\liminf_{k\to\infty }C_{\varepsilon_{n_k}} \exp (-\overline u_{2,\varepsilon_{n_k}}(y_0))
\ge \int_{D}\exp (-u(x,y_0))dx.
\end{equation}
Notice that, from %(\ref{431}), (\ref{4320})  
(\ref{4170})-(\ref{417}) and Lemma 3.7, the following holds:
\begin{align*}
&\int_{D}\exp (-u(x,y_0))dx\\
>&
\exp \left(-\int_{\mathbb{R}^d}\frac{1}{2}|x|^2P_1(dx)+\log \{{\rm Vol}(B_r)\}-\frac{1}{2} \int_{B_r}|x|^2\frac{dx}{{\rm Vol} (B_r)}\right)>0.
\end{align*}
We prove (\ref{4320}).
For sufficiently large $m\ge 1$, 
\begin{align*}
\int_{B\cap D_{R,  \varepsilon }} \exp (-\overline u_{1,\varepsilon}(x)-\overline u_{2,\varepsilon}(y_0))dx
\le &\int_{B} \psi_{\delta, R, \varepsilon }(x)\exp (-\overline u_{1|m,\varepsilon }(x)-\overline u_{2,\varepsilon}(y_0))dx
\end{align*}
(see (\ref{326}) for notation).
Let $\psi_{\delta}$ denote the function $\psi_{\delta, R, 0}$ with $D_{R,0}$ replaced by
$D$.
Then for $m> r$, from (\ref{4170}) and (\ref{4311}),
\begin{align}\label{432}
&\lim_{k\to\infty }\int_{B} \psi_{ \delta, R, \varepsilon_{n_k} }(x)  \exp (-\overline u_{1|m,\varepsilon_{n_k}}(x)-\overline u_{2,\varepsilon_{n_k}}(y_0))dx\\
=&\int_{B}\psi_{\delta, R, 0}(x)\exp (-u_m(x,y_0))dx
\to \int_{B}\psi_{\delta}(x)\exp (-u(x,y_0))dx,\quad m, R\to \infty,\notag\\
\to & \int_{B}I_{D}(x)\exp (-u(x,y_0))dx,\quad \delta\to 0,\notag
\end{align}
since $R\mapsto D_{R,  \varepsilon }$ is nondecrerasing.

We prove (\ref{43201}).
$$\tilde D_{\delta,m,  \varepsilon }:=\{x\in B_r| \overline u_{1,\varepsilon}(x)-\overline u_{1|m,\varepsilon}(x)<\delta\}, \quad \delta>0.$$
Then 
\begin{align}\label{434}
&C_\varepsilon \exp (-\overline u_{2,\varepsilon}(y_0))\\
\ge& \int_{B_r} \psi_{\delta, R,\varepsilon }(x) \exp (-\overline u_{1,\varepsilon}(x)-\overline u_{2,\varepsilon}(y_0))dx\notag\\
\ge &\exp (-\delta )\int_{\tilde D_{\delta,m,  \varepsilon }} \psi_{ \delta, R, \varepsilon }(x) 
\exp (-\overline u_{1|m,\varepsilon}(x)-\overline u_{2,\varepsilon}(y_0))dx\notag\\
=&\exp (-\delta )\int_{B_r} \psi_{\delta,R, \varepsilon }(x)
 \exp (-\overline u_{1|m,\varepsilon}(x)-\overline u_{2,\varepsilon}(y_0))dx\notag\\
&\qquad -\exp (-\delta )\int_{B_r\cap \tilde D_{\delta,m,  \varepsilon }^c} \psi_{\delta, R, \varepsilon}(x) 
\exp (-\overline u_{1|m,\varepsilon }(x)-\overline u_{2,\varepsilon}(y_0))dx.\notag
\end{align}
From (\ref{432}), we only have to prove that the following holds:
\begin{align}\label{435}
\lim_{\delta\to 0}\lim_{R\to\infty}\limsup_{m\to\infty}\limsup_{k\to \infty}
\int_{B_r\cap\tilde D_{\delta,m,  \varepsilon_{n_k} }^c} \psi_{\delta, R,\varepsilon_{n_k}}(x) &\\
\times \exp (-\overline u_{1|m,\varepsilon_{n_k}}(x)-\overline u_{2,\varepsilon_{n_k}}(y_0))dx&=0.\notag
\end{align}
\begin{align*}%\label{436}
&\int_{B_r\cap\tilde D_{\delta,m,  \varepsilon}^c} \psi_{\delta, R,\varepsilon}(x) 
\exp (-\overline u_{1|m,\varepsilon}(x)-\overline u_{2,\varepsilon}(y_0))dx\notag\\
\le &\int_{B_r\cap U_\delta (D_{R,  \varepsilon })\cap D_{R,  \varepsilon }^c}
\exp (-\overline u_{1|m,\varepsilon}(x)-\overline u_{2,\varepsilon}(y_0))dx\notag\\
&\qquad 
+\int_{B_r\cap\tilde D_{\delta,m,   \varepsilon }^c\cap D_{R,  \varepsilon }} 
\exp (-\overline u_{1|m,\varepsilon}(x)-\overline u_{2,\varepsilon}(y_0))dx,\notag
\end{align*}
since $\psi_{\delta, R, \varepsilon}^{-1} ((0,1])=U_\delta (D_{R,  \varepsilon })$.
For any $\gamma>0$, sufficiently large $m\ge m_0\ge 1$ and $k$, from Lemma 3.8 and (\ref{4170}), 
\begin{align*}%\label{437}
&\int_{B_r\cap U_\delta (D_{R,  \varepsilon_{n_k} })\cap D_{R,  \varepsilon_{n_k}  }^c}
\exp (-\overline u_{1|m,\varepsilon_{n_k}}(x)-\overline u_{2,\varepsilon_{n_k} }(y_0))dx\\
\le & \int_{B_r\cap U_{\delta+\gamma } (D_{R,0})\cap U_{-\gamma } (D_{R, 0 })^c}
\exp (-\overline u_{1|m_0,\varepsilon_{n_k}}(x)-\overline u_{2,\varepsilon_{n_k} }(y_0))dx\notag\\
\to & \int_{B_r\cap U_{\delta+\gamma } (D_{R,0})\cap U_{-\gamma } (D_{R, 0 })^c}
\exp (-u_{m_0}(x,y_0))dx,\quad k\to\infty,\notag\\
\to & \int_{B_r\cap U_{\delta+\gamma } (D)\cap U_{-\gamma } (D)^c}\exp (-u_{m_0}(x,y_0))dx,\quad R\to\infty,\notag\\
\to & 0, \quad \delta,\gamma\to 0.\notag
\end{align*}
\begin{align}\label{438}
&\int_{B_r\cap\tilde D_{\delta,m,  \varepsilon }^c\cap D_{R,  \varepsilon }} 
\exp (-\overline u_{1|m,\varepsilon}(x)-\overline u_{2,\varepsilon}(y_0))dx\\
\le&\int_{\tilde D_{\delta,m,  \varepsilon }^c} 
\exp (-\overline u_{1|m,\varepsilon}(x)-\overline u_{2,\varepsilon}(y_0)+R-\overline u_{2,\varepsilon}(y_0))C_\varepsilon p_{0,r,\varepsilon}(x)dx\notag\\
\le & \exp (-2\inf\{\overline u_{1 |m_0,\varepsilon}(x)+\overline u_{2,\varepsilon}(y_0)|x\in B_r \}+R)
{\rm Vol}(B_r)\int_{\tilde D_{\delta,m,  \varepsilon }^c} p_{0,r,\varepsilon}(x)dx,\notag
\end{align}
\begin{align}
&\int_{\tilde D_{\delta,m,  \varepsilon }^c} p_{0,r,\varepsilon}(x)dx=
\int_{\tilde D_{\delta,m,  \varepsilon }^c\times \mathbb{R}^d} \mu_\varepsilon(dxdy)
\le \mu_{1|m,\varepsilon} (\tilde D_{\delta,m,  \varepsilon }^c)+P_{1}(B_m^c)\\
\le & \int_{\tilde D_{\delta,m,  \varepsilon }^c}\exp \left(\frac{\overline u_{1,\varepsilon}(x)-\overline u_{1|m,\varepsilon }(x)-\delta }{\varepsilon}\right)\mu_{1|m,\varepsilon } (dx)+P_{1}(B_m^c)\notag\\
\le&\exp \left(-\frac{\delta }{\varepsilon}\right) \int_{B_r}p_{0,r,\varepsilon}(x)dx+P_{1}(B_m^c)
\to P_{1}(B_m^c),\quad \varepsilon\to 0,\notag\\
\to &0,\quad m\to\infty.\notag
\end{align}
Here, from (\ref{326}) and (\ref{4241}) (see also (\ref{11})), 
\begin{align*}
\exp \left(\frac{\overline u_{1,\varepsilon}(x)-\overline u_{1|m,\varepsilon }(x)}{\varepsilon}\right)\mu_{1|m,\varepsilon } (dx)
=&\exp (u_1 (x;g_\varepsilon (1),  P_{0,r,\varepsilon}, P_1))\nu_{1,\varepsilon} (dx)=p_{0,r,\varepsilon}(x)dx.
\end{align*}
(\ref{416}) and (\ref{4170}) complete the proof of (\ref{435}).

If $P_1$ is compactly supported, then $u_{i|m,\varepsilon}=u_{i|m',\varepsilon}$ and $u(x,y)=u_{m'}(x,y)$
for $m'\ge m$, provided $B_r\cup supp (P_1)\subset B_m$.
(\ref{415})-(\ref{4170}) imply that the last statement of Theorem 2.2 holds.$\Box$

\bigskip
\noindent
Acknowledgement: 
This work was  supported by JSPS KAKENHI Grant Numbers JP26400136 and JP16H03948.
We would also like to thank an anonymous referee for useful suggestions.


\begin{thebibliography}{9}

\bibitem{ADPZ}
S. Adams, N. Dirr, M. A. Peletier, J. Zimmer, 
From a large-deviations principle to the Wasserstein gradient flow: a new micro-macro passage, 
Commun. in Math. Phys. \textbf{307}, (2011), 791--815. 

\bibitem{BCGL}
J. Backhoff, G. Conforti, I. Gentil, C.  L\'eonard,
The mean field Schr\"odinger problem: ergodic behavior, entropy estimates and functional inequalities,
arXiv:1905.02393v1.

\bibitem{Ba}
I. J. Bakelman,
Convex Analysis and Nonlinear Geometric Elliptic Equations,
Springer-Verlag, %Berlin-Heidelberg-New York, 
1994.

\bibitem{S. B}
S. Bernstein, 
Sur les liaisons entre les grandeurs al\'etoires,
Verh. des intern. Mathematikerkongr. Zurich 1932, Band \textbf{1}, (1932), 288--309.

\bibitem{B}
A. Beurling,
An Automorphism of Product Measures,
Ann. of Math. \textbf{72}, (1960), 189--200.

\bibitem{Bi}
P. Billingsley,
Convergence of Probability Measures,
Wiley-Interscience, 1999.

\bibitem{Bre1}
Y. Brenier, 
D\'ecomposition polaire et r\'earrangement monotone des champs de vecteurs, 
C. R. Acad. Sci. Paris S\'erie I, {\bf 305},  no. 19, (1987), 805--808. 

\bibitem{Bre2}
 Y. Brenier,
Polar factorization and monotone rearrangement of vector-valued functions,
Comm. Pure Appl. Math.,  {\bf 44},  no. 4, (1991), 375--417.

\bibitem{C}
Y. Chen, T. Georgiou, M. Pavon,
Entropic and displacement interpolation: a computational approach using the Hilbert metric,
SIAM J. Appl. Math. \textbf{76}, (2016), 2375--2396.

\bibitem{CDS}
M. Colombo, S. D. Marino, F. Stra,
Continuity of multimarginal optimal transport with repulsive cost,
SIAM J. MATH. ANAL.  \textbf{51}, (2019), 2903--2926.

\bibitem{Con0}
G. Conforti, 
A second order equation for Schr\"odinger bridges with applications to the hot gas experiment and entropic transportation cost, 
Probability Theory and Related Fields,  \textbf{174}, (2019), 1--47.

\bibitem{Con}
G. Conforti, L. Ripani, 
Around the entropic Talagrand inequality, 
Bernoulli 26 (2020), no. 2, 1431--1452.
　
\bibitem{Cor}
D. Cordero-Erausquin, B. Klartag,
Moment measures, 
J.Funct. Anal. \textbf{268} (12), (2015), 3834--3866.

\bibitem{7} 
P. Dai Pra,
A stochastic control approach to reciprocal diffusion processes,
Appl. Math. Optim. \textbf{23}, (1991), 313--329.

\bibitem{Doob} 
J.L. Doob,
Conditional Brownian motion and the boundary limits of harmonic functions,
Bulletin de la Soci\'et\'e Math\'ematique de France \textbf{85}, (1957), 431--458.

\bibitem{DLR} 
M. H. Duong, V. Laschos, M. Renger,
Wasserstein gradient flows from large deviations of many-particle limits,
ESAIM: Control, Optimisation and Calculus of Variations, \textbf{19} (4), (2013), 1166--1188.

\bibitem{DE}
P. Dupuis, R. S. Ellis,
A Weak Convergence Approach to the Theory of Large Deviations,
John Wiley \& Sons, %New York, 
1997.

\bibitem{EMR}
M. Erbar,  J. Maas, M. Renger,
From large deviations to Wasserstein gradient flows in multiple dimensions,
Electron. Commun. Probab. \textbf{20}, No. 89, (2015), 1--12.

\bibitem{Jin} 
J. Feng, T. Nguyen,
Hamilton-Jacobi equations in space of measures associated with a system of conservation laws,
J. Math. Pures Appl. \textbf{97}, (2012), 318--390.

\bibitem{F}
W.~H. Fleming,  H.~M. Soner,
Controlled Markov Processes and Viscosity Solutions, 2nd ed.,
Springer, %Berlin-Heidelberg-New York, 
2006.

\bibitem{11} 
H. F\"ollmer,
Random fields and diffusion processes,
in: Hennequin, P.~L.(ed) 
\'Ecole d'\'Et\'e de Probabilit\'es de Saint-Flour XV--XVII, 1985--87,  Lecture Notes in Math.\textbf{1362}, 101--203, 
Springer, %Berlin-Heidelberg-New York, 
2006.
  
\bibitem{Fo} 
R. Fortet,
R\'esolution d'un Syst\`eme d'\'equations de M. Schroedinger,
J. Math. Pures Appl. \textbf{IX}, (1940), 83--105.

\bibitem{GLRT} 
I. Gentila, C. L\'eonard, L. Ripania, L. Tamaninic,
An entropic interpolation proof of the HWI inequality,
Stochastic Process. Appl. \textbf{130}, (2020), 907-923.

\bibitem{J1}
B. Jamison,
Reciprocal Processes,
Z. Wahr. verw. Gebiete \textbf{30} (1974), 65--86.

\bibitem{J2} 
B. Jamison,
The Markov process of Schr\"odinger,
Z. Wahr. Verw. Gebiete \textbf{32} (1975), 323--331.

\bibitem{Leo1}
C. L\'eonard,
From the Schr\"odinger problem to the Monge-Kantorovich problem,
J. Funct. Anal. \textbf{262}, (2012), 1879--1920.

\bibitem{Leo2} 
C. L\'eonard ,
A survey of the Schr\"odinger problem and some of its connections with optimal transport,
Special Issue on Optimal Transport and Applications,
Discrete Contin. Dyn. Syst. \textbf{34}, (2014), 1533--1574.

\bibitem{19}  
T. Mikami,
Monge's problem with a quadratic cost by the zero-noise limit of $h$-path  processes,
Probab. Theory Related Fields \textbf{129}, (2004), 245--260.

\bibitem{21}
T. Mikami,
Semimartingales from  the Fokker-Planck equation.
Appl. Math. Optim. \textbf{53}, (2006), 209--219.

\bibitem{21-1}
T. Mikami:
Marginal problem for semimartingales via duality,
in: Giga, Y. , Ishii, K., Koike, S. et al. (eds)
International Conference for the 25th Anniversary of Viscosity Solutions, 
Gakuto International Series. Mathematical Sciences and Applications \textbf{30}, 133--152,
Gakkotosho, %Tokyo, 
2008.

\bibitem{22} 
T. Mikami,
Optimal transportation problem as stochastic mechanics,
in Selected Papers on Probability and Statistics, Amer. Math. Soc. Transl. Ser. 2,  \textbf{227}, 75--94, 
Amer. Math. Soc., %Providence, RI., 
2009.

\bibitem{23}
T. Mikami,
Two end points marginal problem by stochastic optimal transportation,
SIAM J. Control Optim. \textbf{53}, (2015), 2449--2461.

\bibitem{24}
T. Mikami,
Regularity of Schr\"odinger's functional equation and mean field PDEs for h-path processes,
Osaka J. Math.  \textbf{56}, (2019) 831-842.

\bibitem{25}
T. Mikami, M. Thieullen,
Duality theorem for the stochastic optimal control problem,
Stochastic Process. Appl. \textbf{116}, (2006), 1815--1835.
 
\bibitem{PAL}
S. Pal,
On the difference between entropic cost and the optimal transport cost,
arXiv:1905.12206v1.

\bibitem{RR}
S.T. Rachev,  L. R\"uschendorf, 
Mass transportation problems, Vol. I: Theory, Vol. II: Application, Springer-Verlag, 1998.

\bibitem{R}
L. Ripani,
Convexity and regularity properties for entropic interpolations,
J. Func. Anal. \textbf{277}, (2019), 368--391.

\bibitem{28}
L. R\"uschendorf,  W. Thomsen,
Note on the Schr\"odinger equation and $I$-projections,
Statist. Probab. Lett. \textbf{17}, (1993), 369--375.

\bibitem{Sa}
F. Santambrogio,
Dealing with moment measures via entropy and optimal transport,
J.Funct. Anal. \textbf{271} (2), (2016), 418--436. 

\bibitem{S1}
E. Schr\"odinger,
Ueber die Umkehrung der Naturgesetze,
Sitz. Ber. der Preuss. Akad. Wissen., Berlin, Phys. Math., (1931), p. 144-153.

\bibitem{S2}
E. Schr\"odinger,
Th\'eorie relativiste de l'electron et l'interpr\'etation de la m\'ecanique quantique,
Ann. Inst. H. Poincar\'e \textbf{2}, (1932), 269--310. 

\bibitem{TT}
X. Tan, N. Touzi,
Optimal transportation under controlled stochastic dynamics,
Ann. Probab. \textbf{41}, (2013), 3201--3240.

\bibitem{V}Villani, C.:
\textit{Optimal Transport: Old and New}.Springer-Verlag, Berlin, 2008. 

\bibitem{Z}
J.~C. Zambrini,
Variational processes,
in: Albeverio, S., etal. (eds.)
Stochastic processes in classical and quantum systems (Ascona, 1985), Lecture Notes in Phys. \textbf{262}, 517--529, Springer, %Berlin-Heidelberg-New York, 
1986.


\end{thebibliography}
\end{document}